% !TEX spellcheck = English
\magnification=1000
\def\arXiv{oui}
  % MISE EN PAGE 
  
 \def\oui{oui} 
  
  %\magnification 1100
\ifx\arXiv\oui
\else
 \pdfpagewidth=210truemm
 \pdfpageheight=297truemm 
\fi
  
  %
  % Definitions pour les notes en bas de page
  % [Format : \note{texte de la note} avec numerotation automatique]
  %

  %
  %\anote=\footnote de plain (2 arguments)
  %

  \catcode`@=12 % at signs are no longer letters

  %reference de note
 \def\defrefnote#1{\definexref{#1}{{\the\footnotenumber}}{refnotes}}

  %
  % Fin notes
  %

  %\def\eqlgno#1{$$\leqalignno{#1}$$}

\ifx\couleurs\oui
\input graphicx
 \pdfpagewidth=210truemm
 \pdfpageheight=297truemm 
 \voffset=-5mm
\fi

\input eplain.tex
\expandafter\def\expandafter\newdimen\expandafter{\newdimen}

% couleurs
% Attention : \illust et graphicx nŽcessitent \couleurs={oui}
\ifx\couleurs\oui
%\expandafter\def\expandafter\newdimen\expandafter{\newdimen}
\beginpackages
\usepackage{color}
\endpackages 
 \pdfpagewidth=210truemm
 \pdfpageheight=297truemm 
\long\def\rge#1{{\color{red}#1}}

\definecolor{bleu-iecn}{cmyk}{.98,.13,.1,.55}

\else
\long\def\rge#1{#1}

\fi

\makeatletter
\def\numberedfootnote{%
ÊÊ\global\advance\footnotenumber by 1
ÊÊ\@eplainfootnote{{\number\footnotenumber}}%
}%
\def\makecolumns#1/#2 {\par \begingroup
ÊÊ \@columndepth = #1
ÊÊ \advance\@columndepth by -1
ÊÊ \divide \@columndepth by #2
ÊÊ \advance\@columndepth by 1
ÊÊ \@linestogoincolumn = \@columndepth
ÊÊ \@linestogo = #1
ÊÊ \currentcolumn = 1
ÊÊ \def\@endcolumnactions{%
ÊÊÊÊÊÊ\ifnum \@linestogo<2
ÊÊÊÊÊÊÊÊ \the\crtok \egroup \endgroup \par % End \valign and \makecolumns.
ÊÊÊÊÊÊ\else
ÊÊÊÊÊÊÊÊ \global\advance\@linestogo by -1
ÊÊÊÊÊÊÊÊ \ifnum\@linestogoincolumn<2
ÊÊÊÊÊÊÊÊÊÊÊÊ\global\advance\currentcolumn by 1
ÊÊÊÊÊÊÊÊÊÊÊÊ\global\@linestogoincolumn = \@columndepth
ÊÊÊÊÊÊÊÊÊÊÊÊ\the\crtok
ÊÊÊÊÊÊÊÊ \else
ÊÊÊÊÊÊÊÊÊÊÊÊ&\global\advance\@linestogoincolumn by -1
ÊÊÊÊÊÊÊÊ \fi
ÊÊÊÊÊÊ\fi
ÊÊ }%
ÊÊ \makeactive\^^M
ÊÊ \letreturn \@endcolumnactions
ÊÊ \@columnwidth = \hsize
ÊÊÊÊ \advance\@columnwidth by -\parindent
ÊÊÊÊ \divide\@columnwidth by #2
ÊÊ \penalty\abovecolumnspenalty
ÊÊ \noindent % It's not a paragraph (usually).
ÊÊ \valign\bgroup
ÊÊÊÊ &\hbox to \@columnwidth{\strut \hsize = \@columnwidth ##\hfil}\cr
}%
\makeatother

\lefteqnumbers
% choix numerotation gauche ou droite: defaut=gauche, \numadroite{oui} -> a droite
   \def\testd{oui}
   \def\choixlat{\ifx\numadroite\testd\righteqnumbers
            \else  \lefteqnumbers\fi}
    \choixlat

% modif format note
\catcode`@=\letter
\def\@eplainfootnote#1{\let\@sf\empty % parameter #2 (the text) is read later
  \ifhmode\edef\@sf{\spacefactor\the\spacefactor}\/\fi
  \global\advance\hlfootlabelnumber by 1
  \hlstart@impl{foot}{\hlfootlabel}%
  \hldest@impl{footback}{\hlfootbacklabel}%
  \hbox{$^{(#1)}$}%
  \hlend@impl{foot}%
  \@sf\vfootnote{#1.}%
}%
\catcode`@=\other

  \interfootnoteskip=0pt
  
  \everyfootnote={\eightpoint\leftskip=5truemm\rightskip5truemm}
  
  \hsize150truemm\vsize 240truemm\hoffset=5truemm

  \pretolerance=500\tolerance=1000\brokenpenalty=5000
  \parindent3mm
  
  \countdef\temps=170
  \temps=\time
  \countdef\nminutes=171{\nminutes = \time}
  \countdef\nheure=172
  \def\heure{\begingroup                     % heure a la francaise
     \temps = \time \divide\temps by 60
     \nheure = \temps                        % l'heure, de 0 \`a 23
     \nminutes = \time
     \multiply\temps by 60
     \advance\nminutes by -\temps            % Les minutes, de 0 a 59.
     \ifnum\nminutes<10 \toks1 = {0}%
     \else\toks1 = {}%
     \fi
     \number\nheure h\the\toks1 \number\nminutes  
  \endgroup}%

  \newcount\chstart
  \chstart=\pageno
 \headline={\ifnum\pageno=\chstart {\hfill} \else {\hss \tenrm --\ \folio\ --\hss}\fi}
  \footline={\hfill}
  \normalbaselines
  \frenchspacing
    \def\dater{\vglue-10mm\rightline{(\the\day/\the\month/\the\year)}}
  \def\dateheure{\vglue-10mm\rightline{(\the\day/\the\month/\the\year,\ \heure)}}
 % justification ˆ droite

  \newif\ifpagetitre \pagetitretrue
\newtoks\hautpagetitre \hautpagetitre={\hfill}
\newtoks\baspagetitre \baspagetitre={\hfill}
\newtoks\auteurcourant \auteurcourant={\hfill}
\newtoks\titrecourant \titrecourant={\hfill}
\newtoks\hautpagegauche
\newtoks\hautpagedroite
\newtoks\hautpagemilieu
\hautpagemilieu={\tenrm\hfil -- \folio\ -- \hfil}
\hautpagegauche={\ifx\midfolio\oui\the\hautpagemilieu\else\tenrm\folio\hfill\the\auteurcourant\hfill\fi}
\hautpagedroite={\ifx\midfolio\oui\the\hautpagemilieu\else\hfill\the\titrecourant\hfill\tenrm\folio\fi}
\newtoks\baspagegauche \baspagegauche={\hfil}
\newtoks\baspagedroite \baspagedroite={\hfil}
\headline={\ifpagetitre\the\hautpagetitre
\else\ifodd\pageno\the\hautpagedroite\else\the\hautpagegauche\fi\fi }
\footline={\ifpagetitre\the\baspagetitre
\else\ifodd\pageno\the\baspagedroite
\else\the\baspagegauche\fi\fi \global\pagetitrefalse}

\def\pageblanche{\vfill\eject\pagetitretrue
\null\vfill\eject
\pagetitretrue
}
\def\chgtpage{\ifodd\pageno \else
\pageblanche \fi \pagetitretrue\titreun=0\footnotenumber=0}

\def\chgtpageincrtitreun{\ifodd\pageno \else
\pageblanche \fi \pagetitretrue\footnotenumber=0}

\def\majnombres{\ifodd\pageno \else
\pageblanche \fi \pagetitretrue\hautpoly\titreun=0\footnotenumber=0}

\def\hautspages#1#2{\auteurcourant={\ninepcap#1}\titrecourant={\nineit#2}}

\ifnum\chstart=\pageno \pagetitretrue\fi
%\ifnum\pageno=\chstartÊ\pagetitretrue \fi
  
  %pour eviter les rectangles noirs:
  %\overfullrule=Opt
  % pour que les @ soient des lettres

  %% pour aller a la ligne dans un \proclaim:
  \def\PAR{\par}
  
  %% Note en marge gauche
  \def\leftnote#1{\vadjust{\setbox1=\vtop{\hsize 20mm\parindent=0pt\eightpoint
  \baselineskip=9pt\rightskip=4mm plus 4mm\vskip-4.7mm#1}\hbox{\kern-2cm\smash{\box1}}}}

  %% Pour couper automatiquement un titre trop long
  
  \def\raggedcenter{\leftskip=20pt plus 10em  % reglages initial 4em
       \rightskip=\leftskip 
        \parfillskip=0pt 
         \spaceskip=.3333em \xspaceskip=.5em 
          \pretolerance=9999 \tolerance=9999
           \hyphenpenalty=9999 \exhyphenpenalty=9999 }
           
  \def\titrecentre#1{{\parindent0mm\raggedcenter
       \spaceskip=.6em plus .2em minus .2em\xspaceskip=.6em plus .2em minus .2em
        \tit#1\par}}
        
% entourer un texte par un filet

  %%%%%% D\'ebut num\'erotation automatique des paragraphes et/ou des enonces

  \def\oui{oui}
  
\def\fontetitreun{\ifx\paradouze\oui\douzepts\gpdouze\twelvebf\textfont1=\twelveib\else
\quatorzepts\gpquatorze\fourteenbf\fi}

\def\fontetitreunl{\douzepts\textfont1=\twelveib\scriptfont1=\tenib\fourteenti}
 
 \def\fontetitredeux{\textfont1=\eleveni\ifx\paradouze\oui\onzepts\scriptfont1=\ninei\elevenit\else
                        \douzepts\twelveit\fi}
 
   \def\fontetitredeuxb{\ifx\paradouze\oui\onzepts\eleventi\gponze\textfont1=\elevenib\scriptfont1=\nineib
                         \else\douzepts\twelveti\scriptfont1=\twelveib\scriptfont1=\tenib\gpdouze\fi}
                         
\def\fontetitredeuxl{\onzepts\textfont1=\elevenbf\scriptfont1=\ninebf\twelvebf}
  
\def\fontetitretrois{\textfont0=\elevenrm\scriptfont0=\eightrm\textfont1=\eleveni
                      \scriptfont1=\eighti\scriptscriptfont1=\sixi\elevenit}
                      
\def\fontetitrequatre{\textfont0=\elevenrm\scriptfont0=\eightrm\textfont1=\eleveni
                      \scriptfont1=\eighti\scriptscriptfont1=\sixi\elevenrm}
  
  \newcount\titreun\titreun=0
  \newcount\titredeux\titredeux=0
  \newcount\titretrois\titretrois=0
  \newcount\titrequatre\titrequatre=0
  \newcount\enonce\enonce=0
  
  \def\incr#1{\global\advance#1 by 1 {\the #1}}
  \def\avance#1{\global\advance#1 by 1}
  \def\init#1{\global#1=0}
  
  \long\def\Indentation#1#2{\setbox10=\hbox{\fontetitreun#1}
  	                    \ifdim\wd10 < 4mm
                         \setbox10=\hbox to 4mm{\box10\hfill}
                       \else\ifdim\wd10 < 6mm
                         \setbox10=\hbox to 6mm{\box10\hfill}
  	                    \else\ifdim\wd10 < 8mm
                         \setbox10=\hbox to 8mm{\box10\hfill}
                       \else\ifdim\wd10 < 12mm
                         \setbox10=\hbox to 12mm{\box10\hfill}
                       \fi\fi\fi\fi
                       \dimen10=\hsize
                       \advance \dimen10 by -\wd10
                       \noindent \box10 %
                       \ignorespaces
                       \hbox{\vtop{\hsize=\dimen10\raggedright\noindent\fontetitreun#2}}}

  \long\def\paraun#1{\removelastskip\par\medskip\goodbreak\vskip0pt plus.01\vsize\penalty-100
                \vskip0pt plus-.01\vsize
  	              \init{\titredeux}\ifnum\optionparag=1{\init\eqnumber\init\enonce}\else{}\fi
                  \goodbreak{\fontetitreun
  	                \Indentation{\incr{\titreun}.\ }{\fontetitreun #1\par}}\nobreak\medskip}

 %
 % titres de paragraphes centres
 %
 \long\def\paraunc#1{\removelastskip\par\bigskip\goodbreak\vskip0pt plus.01\vsize\penalty-100
                \vskip0pt plus-.01\vsize
  	              \init{\titredeux}
                 \ifnum\optionparag=1{\init{\eqnumber}\init\enonce}\else{}\fi
                  \goodbreak
  	                {\parindent0mm\raggedcenter\fontetitreun\incr{\titreun}.\ 
                     \fontetitreun #1\par}\nobreak\medskip}
                     
% paragraphes de niveau 1 numŽrotŽs avec des chiffres romains
\newtoks\titreunl
\titreunl={\ifnum\titreun=1{I}\fi%
\ifnum\titreun=2{II}\fi%
\ifnum\titreun=3{III}\fi%
\ifnum\titreun=4{IV}\fi%
\ifnum\titreun=5{V}\fi%
\ifnum\titreun=6{VI}\fi%
\ifnum\titreun=7{VII}\fi%
\ifnum\titreun=8{VIII}\fi%
\ifnum\titreun=9{IX}\fi%
\ifnum\titreun=10{X}\fi%
\ifnum\titreun=11{XI}\fi%
\ifnum\titreun=12{XII}\fi%
\ifnum\titreun=13{XIII}\fi%
}
\long\def\paraunl#1{\removelastskip\par\bigskip\bigskip\goodbreak\vskip0pt plus.01\vsize\penalty-100
                \vskip0pt plus-.01\vsize
  	              \init{\titredeux}\ifnum\optionparag=1{\init\eqnumber\init\enonce}\else{}\fi
                  \goodbreak{\fontetitreunl
  	                \Indentation{\global\advance\titreun by 1{\the\titreunl}.\ }{\fontetitreunl #1\par}}\nobreak\smallskip}

% paragraphes de niveaux 2
  
  \long\def\paradeux#1{\init{\titretrois}\vskip0pt plus.01\vsize\penalty-10
                \vskip0pt plus-.01\vsize\ifx \elie\oui\medskip\ifnum\titredeux>0\medskip\fi\fi
                 \Indentation{\fontetitredeux\the\titreun${\cdot}$\incr{\titredeux}.}
                              {\fontetitredeux\textfont1=\eleveni#1}\nobreak\par }
  
  \long\def\paradeuxb#1{\init{\titretrois}\vskip0pt plus.001\vsize\penalty-10
                \vskip0pt plus-.01\vsize{\ifx \elie\oui\medskip\ifnum\titredeux>0\medskip\fi\fi
                  \Indentation
  {\fontetitredeuxb\the\titreun${\cdot}$\incr{\titredeux}.}{ \fontetitredeuxb#1}}\nobreak
\smallskip}

% paragraphes de niveau 2 numŽrotŽs avec des lettres capitales 
\newtoks\titredeuxl
\titredeuxl={\ifnum\titredeux=1{A}\fi%
\ifnum\titredeux=2{B}\fi%
\ifnum\titredeux=3{C}\fi%
\ifnum\titredeux=4{D}\fi%
\ifnum\titredeux=5{E}\fi%
\ifnum\titredeux=6{F}\fi%
\ifnum\titredeux=7{G}\fi%
\ifnum\titredeux=8{H}\fi%
\ifnum\titredeux=9{I}\fi%
\ifnum\titredeux=10{J}\fi%
\ifnum\titredeux=11{K}\fi%
\ifnum\titredeux=12{L}\fi%
\ifnum\titredeux=13{M}\fi%
}
 \long\def\paradeuxl#1{\init{\titretrois}\vskip0pt plus.001\vsize\penalty-10
                \vskip0pt plus-.01
                \vsize \bigskip%
                  \Indentation
     {\fontetitredeuxl\global\advance\titredeux by 1
  \quad \the\titreunl${\cdot}$\the\titredeuxl.}{ \fontetitredeuxl#1}
  \removelastskip\nobreak\smallskip}
  
% paragraphes de niveaux 3

  \long\def\paratrois#1{\init{\titrequatre}\ifdim\lastskip<\smallskipamount
                \removelastskip\smallskip\fi
                 \vskip0pt plus.01\vsize\penalty-10
                  \vskip0pt
plus-.01\vsize{\ifx \elie\oui\ifnum\titretrois>0\medskip\fi\fi
\Indentation{\fontetitretrois\the\titreun${\cdot}$\the\titredeux${\cdot}$\incr{\titretrois}.\ }
  {\hskip0mm\baselineskip=14pt\fontetitretrois#1}\nobreak\smallskip}}
  
 % paragraphes de niveau 3 numŽrotŽs avec des chiffres mais niveaux 1 et 2 chiffres romains et lettres 
  
  \long\def\paratroisl#1{\init{\titrequatre}\ifdim\lastskip<\smallskipamount
                \removelastskip\fi
                 \vskip0pt plus.01\vsize\penalty-10
                  \vskip0pt
plus-.01\vsize\ifx \elie\oui\bigskip
\fi
\Indentation{\fontetitretrois\quad \quad \the\titreunl{${\cdot}$}\the\titredeuxl${\cdot}$\incr{\titretrois}.\ }
  {\hskip0mm\fontetitretrois#1}\nobreak\smallskip}

% paragraphes de niveaux 4

  \long\def\paraquatre#1{\ifdim\lastskip<\smallskipamount
                \removelastskip\smallskip\fi
                 \vskip0pt plus.01\vsize\penalty-10
                  \vskip0pt
                  plus-.01\vsize\par
                %\bigskip
 
\Indentation{\fontetitrequatre \the\titreun{${\cdot}$}\the\titredeux${\cdot}$\the\titretrois${\cdot}$\incr{\titrequatre}.\ }
{\hskip0mm\fontetitrequatre#1}\nobreak\smallskip}

% paragraphes de niveau 4 numŽrotŽs avec des lettres minuscules

\newtoks\titrequatrel
\titrequatrel={\ifnum\titrequatre=1{a}\fi%
\ifnum\titrequatre=2{b}\fi%
\ifnum\titrequatre=3{c}\fi%
\ifnum\titrequatre=4{d}\fi%
\ifnum\titrequatre=5{e}\fi%
\ifnum\titrequatre=6{f}\fi%
\ifnum\titrequatre=7{g}\fi%
\ifnum\titrequatre=8{h}\fi%
\ifnum\titrequatre=9{i}\fi%
\ifnum\titrequatre=10{j}\fi%
\ifnum\titrequatre=11{k}\fi%
\ifnum\titrequatre=12{l}\fi%
\ifnum\titrequatre=13{m}\fi%
}
\long\def\paraquatrel#1{\ifdim\lastskip<\smallskipamount
                \removelastskip\smallskip\fi
                 \vskip0pt plus.01\vsize\penalty-10
                  \vskip0pt
                  plus-.01\vsize{\bigskip
\Indentation{\global\advance\titrequatre by 1
\fontetitrequatre\quad \quad \quad \the\titreunl${\cdot}$\the\titredeuxl${\cdot}$\the\titretrois${\cdot}$\the\titrequatrel.\ }
{\hskip0mm\fontetitrequatre#1}\nobreak\smallskip}}

% Pour memoriser le no de section d'un titre de paragraphe
\ifx\optionkeys\oui
\def\drefun#1{\definexref{¤#1}{{\the\titreun}}{}} 
\def\drefdeux#1{\definexref{¤#1}{{\the\titreun}.{\the\titredeux}}{}}
\def\dreftrois#1{\definexref{¤#1}{{\the\titreun}.{\the\titredeux}.{\the\titretrois}}{}}
\else
\def\drefun#1{\definexref{prg#1}{{\the\titreun}}{}} 
\def\drefdeux#1{\definexref{prg#1}{{\the\titreun}.{\the\titredeux}}{}}
\def\dreftrois#1{\definexref{prg#1}{{\the\titreun}.{\the\titredeux}.{\the\titretrois}}{}}
\fi

% Dans le cas de \drefun, coder \paraun{Titre}\drefun{label} dans le cas de \drefdeux, on peut coder
%                                \paradeux{Sous-titre\drefdeux{label}}
% pour appeler la reference "label", taper \refn{¤label}
%
% en cas de probl\`emes avec option_keys: coder \def\optionkeys{oui} et remplacer ¤ par prg dans \refn

%% num\'erotation des \'enonc\'es
%%Exemples
%%\propt{leb}{ (Lebesgue)}{ L'espace $L^1$ est complet.}
%%\par
%%Ca marche, comme le dit le \ref{leb}.
%%\Propt{leb2}{L'espace $L^1$ est vraiment complet.}\par
%%Ce \ref{leb2} est moins bien que les Th\'eor\`emes \ref{sleb} et \ref{sleb2}.

  \long\def\propdeux#1#2#3#4{%
       \avance{\enonce}
       \leavevmode\edef\temp{#2}%
         \ifx\temp\empty 
          \else
           \definexref{#2}{#1~{\the\titreun.\the\enonce}}{enonces}
            \definexref{s#2}{{\the\titreun.\the\enonce}}{enonces}
             \fi
\smallskip
      \noindent{\bf#1\ {\bf\the\titreun.\the\enonce{#3}.}\enspace}{\sl#4\par}%
      \ifdim\lastskip<\medskipamount \removelastskip\penalty55\par \fi
   }

  \long\def\propun#1#2#3#4{%
      \avance{\enonce}
       \leavevmode\edef\temp{#2}%
        \ifx\temp\empty 
          \else
           \definexref{#2}{#1~{\the\enonce}}{enonces}
            \definexref{{s#2}}{{\the\enonce}}{enonces}
             \fi
   \par 
     \noindent{\bf#1\ {\bf\the\enonce{#3}.}\enspace}{\sl#4\par}%
     \ifdim\lastskip<\medskipamount \removelastskip\penalty55\medskip\fi
  }
  
  \long\def\prop#1#2#3#4{\ifnum\optionparag=1
                          \propdeux{#1}{#2}{\textfont1=\elevenib#3}{#4} \else\propun{#1}{#2}{\textfont1=\elevenib#3}{#4}\fi}

  \long\def\propt#1#2#3{\ifx\tpf\oui \prop{Th\'eo\-r\`eme}{#1}{#2}{#3}\par
                       \else\prop{Theorem}{#1}{#2}{#3}\par\fi}
  \long\def\Propt#1#2{\propt{#1}{}{#2}}
  \long\def\propl#1#2#3{\ifx\tpf\oui\prop{Lem\-me}{#1}{#2}{#3}\par
                         \else\prop{Lemma}{#1}{#2}{#3}\par\fi}
  
  \long\def\propc#1#2#3{\ifx\tpf\oui\prop{Corol\-laire}{#1}{#2}{#3}\par
                         \else\prop{Corollary}{#1}{#2}{#3}\par\fi}
  \long\def\Propc#1#2{\propc{#1}{}{#2}}

  \long\def\propd#1#2#3{\ifx\tpf\oui\prop{D\'efi\-nition}{#1}{#2}{#3}\par
                       \else\prop{Definition}{#1}{#2}{#3}\par\fi} 
  
  \long\def\proptd#1#2#3{\ifx\tpf\oui\prop{Th\'eor\`eme et d\'efi\-nition}{#1}{#2}{#3}\par
                       \else\prop{Theorem and definition}{#1}{#2}{#3}\par\fi}

  % Pour une numerotation des enonces et des formules avec un seul nombre 
  % coder \optionparag=2

  % Choisir numerotation manuelle ou automatique des paragraphes: 
  % utiliser \section au lieu de \paraun et choisir \optionparag=2
  % Cela induit aussi une num\'erotation des formules avec un seul nombre
  
  \newcount\optionparag\optionparag=1
  
  \long\def\section#1#2{\ifnum\optionparag=1 \paraun{#2} 
                        \else\goodbreak{\fontetitreun
  	                \Indentation{#1.\ }{#2}}\nobreak\smallskip\fi}

    % numerotation automatique (par defaut)
    % numerotation manuelle
  \def\eqconstruct#1{\ifnum\optionparag=1{\the\titreun\hbox{$\cdot$}#1}\else{#1}\fi}

  %%%%%%% Fin num\'erotation automatique des paragraphes
  
  %% num\'erotation bibliographie
  
  \def\numref{oui}  % numeroter la biblio par defaut avec bibtem
  
  \newcount\mesref\mesref=0 
  \def\defbib#1{\ifx\numref\oui\global\advance\mesref by 1 \definexref{#1}{{\the
                 \mesref}}{}\else\definexref{#1}{#1}{}\fi}
  \def\bibtem#1{\defbib{#1}\item{\citer{#1}}}
  \def\citer#1{[\ref{#1}]}
  \def\citeplus#1#2{[\ref{#1}; #2]}

  % FONTES & FAMILLES
  
  \font\seventeenmsa=msam10 at 17pt    % symboles d'AMSTEX
  \font\fourteenmsa=msam10 at 14pt
  \font\twelvemsa=msam10 at 12pt
  \font\tenmsa=msam10                 
  \font\ninemsa=msam10 at 9pt 
  \font\eightmsa=msam10 at 8pt 
  \font\sevenmsa=msam7 
  \font\sixmsa=msam10 at 6pt
  \font\fivemsa=msam5
  \newfam\msafam\textfont\msafam=\tenmsa\scriptfont\msafam=\sevenmsa\scriptscriptfont\msafam=\fivemsa
  
  \font\seventeenbb=msbm10 at 17pt     % Lettres evidees pour titres
  \font\fourteenbb=msbm10 at 14pt
  \font\twelvebb=msbm10 at 12pt
  \font\tenbb=msbm10                   %Lettres evidees
  \font\ninebb=msbm10 at 9pt 
  \font\eightbb=msbm10 at 8pt 
  \font\sevenbb=msbm7 
  \font\sixbb=msbm10 at 6pt
  \font\fivebb=msbm5 
  \newfam\bbfam\textfont\bbfam=\tenbb\scriptfont\bbfam=\sevenbb\scriptscriptfont\bbfam=\fivebb
  \def\bb{\fam\bbfam\tenbb}%

  \font\seventeenscaln=eusm10 at 17pt   % Lettres de ronde
  \font\twelvescaln=eusm10 at 12pt
  \font\tenscaln=eusm10                
  \font\ninescaln=eusm10 scaled 900
  \font\eightscaln=eusm10 scaled 800
  \font\sevenscaln=eusm10 scaled 700
  \font\sixscaln=eusm10 scaled 600
   
  \newfam\scalnfam\textfont\scalnfam=\tenscaln\scriptfont\scalnfam=\sevenscaln\scriptscriptfont\scalnfam=\sixscaln
  \def\scaln{\fam\scalnfam\tenscaln}%
  \def\scal{\scaln}
  
  \font\tenscalb=eusb10                % Lettres de ronde grasses

  \font\sevenscalb=eusb10 scaled 700

  \newfam\scalbfam\textfont\scalbfam=\tenscalb\scriptfont\scalbfam=\sevenscalb
  %
  
  %
  % Romain
  %
  \font\fourteenrm=cmr12 scaled 1200
  \font\elevenrm=cmr10 at 11pt
  \font\twelverm=cmr12
  \font\ninerm=cmr9
  \font\eightrm=cmr8      
  \font\sevenrm=cmr7
  \font\sixrm=cmr6

  \font\seventeenpcap=cmcsc10 at 17pt
  \font\tenpcap=cmcsc10                        % Petites capitales
  \font\ninepcap=cmcsc9
  \font\eightpcap=cmcsc8
  \font\sevenpcap=cmcsc10 scaled 700
  
  \newfam\pcapfam\textfont\pcapfam=\tenpcap\scriptfont\pcapfam=\sevenpcap
  \def\pcap{\fam\pcapfam\tenpcap}
  
               % Gras romain (boldface)
                % Titres
  \font\seventeenrm=cmbx12 scaled 1400

  \font\fourteenbf=cmbx10 scaled 1400
  
  \font\twelvebf=cmbx12
  \font\elevenbf=cmbx10 at 11pt
  \font\ninebf=cmbx9  
  \font\eightbf=cmbx8
  \font\sixbf=cmbx6
  
  \font\tengot=eufm10                           % Lettres gothiques
   
  \font\eightgot=eufm10 at 8truept 
  \font\sevengot=eufm7 
  \font\sixgot=eufm10 at 6 truept 
   
  \newfam\gotfam
  \textfont\gotfam=\tengot\scriptfont\gotfam=\sevengot\scriptscriptfont\gotfam=\sixgot
  \def\got{\fam\gotfam\tengot}%

  %% Pour les titres (168pt)
  
  \def\tit{%
  \textfont0=\seventeenrm\scriptfont0=\tenrm\def\rm{\fam0\seventeenrm}%
  \textfont1=\seventeenib\scriptfont1=\twelveib%
  \textfont2=\seventeensy\scriptfont2=\twelvesy\scriptscriptfont2=\ninesy
  \textfont3=\seventeenex
  \textfont\itfam=\seventeenti
  \def\it{\fam\itfam\seventeenti}%
  \textfont\bbfam=\seventeenbb \scriptfont\bbfam=\twelvebb
  \def\bb{\fam\bbfam\seventeenbb}%
  \textfont\msafam=\seventeenmsa\scriptfont\msafam=\twelvemsa
  \textfont\scalnfam=\seventeenscaln
  \def\pcap{\fam\pcapfam\seventeenpcap}
  \normalbaselineskip=25pt\normalbaselines\rm}

  % italiques grasses pour titres
  \font\seventeenti=cmbxti10 scaled 1680
  
  \font\fourteenti=cmbxti10 at 14pt
  
  \font\twelveti=cmbxti10 scaled 1200
  \font\eleventi=cmbxti10 at 11pt

  %
  % italiques
  %
  \font\twelveit=cmti12	
  \font\elevenit=cmti10 scaled 1100
  \font\nineit=cmti9
  \font\eightit=cmti8
  \font\sevenit=cmti7

  %
  % italique mathematique gras pour titres
  %
 
 \font\seventeenib=cmmib10 scaled 1680
  \font\fourteenib=cmmib10 scaled 1400
  \font\twelveib=cmmib10 scaled 1200
  \font\elevenib=cmmib10 scaled 1100
  \font\tenib=cmmib10
\font\eightib=cmmib10 scaled 800
  \font\nineib=cmmib10 scaled 900
\font\sevenib=cmmib10 scaled 700
\font\sixib=cmmib10 scaled 600
\font\fiveib=cmmib10 scaled 500

\ifx\ITAN\oui
\else
\innernewfam\cmmibfam
\textfont\cmmibfam=\tenib
\scriptfont\cmmibfam=\sevenib
\scriptscriptfont\cmmibfam=\fiveib
\def\ib{\fam\cmmibfam\tenib}
\fi

  %
  % italique mathematique
  %
  \font\twelvei=cmmi10 scaled 1200
  \font\eleveni=cmmi10 scaled 1100
  \font\ninei=cmmi9
  \font\eighti=cmmi8 
  \font\seveni=cmmi7 			                
  \font\sixi=cmmi6
  
  \font\ninesl=cmsl9                    % slanted 
  \font\eightsl=cmsl8 
  \font\sevensl=cmsl10 at 7pt

  \font\ninett=cmtt9                    % typewriter
  \font\eighttt=cmtt8
  \font\seventt=cmtt10 scaled 700

  \font\seventeensy=cmsy10 scaled 1680    % symboles pour titres
  \font\fourteensy=cmsy10 scaled 1400
  \font\twelvesy=cmsy10 scaled 1176
  
  \font\ninesy=cmsy9                      % symboles
  \font\eightsy=cmsy8
  \font\sixsy=cmsy6
  \font\seventeenex=cmex10 at 17pt
  \font\fourteenex=cmex10 at 14pt
  \font\twelveex=cmex10 at 12pt
  \font\nineex=cmex10 at 9pt
  \font\eightex=cmex10 at 8pt
  \font\sevenex=cmex10 at 7pt
  \font\sixex=cmex10 at 6pt
  \font\fiveex=cmex10 at 5pt
  
   %% Lettres grecques plombees              %% Lettres grecques plombees
   
  \font\fourteengp=cmmi10 at 14pt
  
  \font\twelvegp=cmmib10 at 12pt
  \font\elevengp=cmmib10 at 11pt
  \font\tengp=cmmib10                          
  \font\ninegp=cmmib10 at 9pt 
  \font\eightgp=cmmib8 
  \font\sevengp=cmmib7 
  \font\sixgp=cmmib6

%\innernewfam\gpfam
%\textfont1\gpfam=\tengp
%\scriptfont1\gpfam=\sevengp
%\scriptscriptfont1\gpfam=\fivegp
%\def\gp{\fam\gpfam\tengp}

  \def\gponze{\textfont0=\elevenbf\scriptfont0=\eightbf\scriptscriptfont0=\sixbf
           \textfont1=\elevengp\scriptfont1=\eightgp\scriptscriptfont1=\sixgp}
  \def\gpdouze{\textfont0=\twelvebf\scriptfont0=\tenbf\scriptscriptfont0=\ninebf
           \textfont1=\twelvegp\scriptfont1=\tengp\scriptscriptfont1=\ninegp}        
  
 \def\gpquatorze{\textfont0=\fourteenbf\scriptfont0=\twelvebf\scriptscriptfont0=\elevenbf
           \textfont1=\fourteengp\scriptfont1=\twelvegp\scriptscriptfont1=\elevengp}

  % CARACTERES SPECIAUX
  
  \expandafter\chardef\csname pre amssym.def at\endcsname=\the\catcode`\@
  \catcode`\@=11
  \def\undefine#1{\let#1\undefined}
  \def\newsymbol#1#2#3#4#5{\let\next@\relax
   \ifnum#2=\@ne\let\next@\msafam@\else
   \ifnum#2=\tw@\let\next@\bbfam@\fi\fi
   \mathchardef#1="#3\next@#4#5}
  \def\mathhexbox@#1#2#3{\relax
   \ifmmode\mathpalette{}{\m@th\mathchar"#1#2#3}%
   \else\leavevmode\hbox{$\m@th\mathchar"#1#2#3$}\fi}
  \def\hexnumber@#1{\ifcase#1 0\or 1\or 2\or 3\or 4\or 5\or 6\or 7\or 8\or
   9\or A\or B\or C\or D\or E\or F\fi}
  
  \def\setboxz@h{\setbox\z@\hbox}
  \def\wdz@{\wd\z@}
  \def\boxz@{\box\z@}
  
  \edef\msafam@{\hexnumber@\msafam}
  \mathchardef\dabar@"0\msafam@39
  
  \edef\bbfam@{\hexnumber@\bbfam}
  \def\widehat#1{\setboxz@h{$\m@th#1$}%
   \ifdim\wdz@>\tw@ em\mathaccent"0\bbfam@5B{#1}%
   \else\mathaccent"0362{#1}\fi}
  \def\widetilde#1{\setboxz@h{$\m@th#1$}%
   \ifdim\wdz@>\tw@ em\mathaccent"0\bbfam@5D{#1}%
   \else\mathaccent"0365{#1}\fi}
  \newsymbol\leqq 1335          % superieur ou egal(=)
  \newsymbol\leqslant 1336
  \newsymbol\lessgtr 1337       % superieur ou inferieur
  \newsymbol\backprime 1038     % apostrophe de gauche a droite
  \newsymbol\risingdotseq 133A  % egal entre points (bas puis haut)
  \newsymbol\fallingdotseq 133B % egal entre points (haut puis bas)
  \newsymbol\succcurlyeq 133C   % superieur ou egal tordu
  \newsymbol\geqq 133D          % inferieur ou egal(=)
  \newsymbol\geqslant 133E
  \newsymbol\nmid 232D
  \newsymbol\nexists 2040
  \newsymbol\smallsetminus 2272
  \newsymbol\varnothing 203F 
  \catcode`\@=\active

  % typographie francaise ou anglaise
  \catcode`\@=11
  \newcount\typofr\typofr=1
  
  \catcode`\;=\active
  \def;{\ifnum\typofr=1\relax\ifhmode\ifdim\lastskip>\z@\unskip\fi
     \kern.2em\fi\string;\else\string;\fi}
  
  \catcode`\:=\active
  \def:{\ifnum\typofr=1\relax\ifhmode\ifdim\lastskip>\z@\unskip\fi
  \penalty\@M\ \fi\string:\else\string:\fi}
  
  \catcode`\!=\active
  \def!{\ifnum\typofr=1\relax\ifhmode\ifdim\lastskip>\z@\unskip\fi
     \kern.2em\fi\string!\else\string!\fi}
  
  \catcode`\?=\active
  \def?{\ifnum\typofr=1\relax\ifhmode\ifdim\lastskip>\z@\unskip\fi
     \kern.2em\fi\string?\else\string?\fi}

  \def\francais{\typofr=1\def\tpf{oui}}
  \def\anglais{\typofr=2\def\tpf{non}\def\english{oui}}
  \def\oui{oui}
  \francais
  
  \catcode`\@=12
  %\catcode`\@=\active
  
  %pour rayer du texte en rouge

\ifx\textures\oui
\def\raye #1|{\leavevmode\setbox1=\hbox{#1}%
\raise .5pt\hbox to \wd1{\xleaders\hbox{\rge{$ \sct / $}%
\kern 1pt}\hfill\kern -1pt }\kern -\wd1 \unhbox1\relax }

\def\barre #1|{\leavevmode\setbox1=\hbox{#1}%
\rlap{\Red\vrule height 2.4pt depth -1.2pt width \wd1}\Black \unhbox1\relax}
\else
\def\raye #1|{\leavevmode\setbox1=\hbox{#1}%
\raise .5pt\hbox to \wd1{\xleaders\hbox{\rge{$ \sct / $}%
\kern 1pt}\hfill\kern -1pt }\kern -\wd1 \unhbox1\relax }

\def\barre #1|{\leavevmode\setbox1=\hbox{#1}%
\rlap{\color{red}\vrule height 2.4pt depth -1.2pt width \wd1}\color{black} \unhbox1\relax}

\fi
  
  % FIN CARACTERES SPECIAUX

  % MACROS DIVERSES
  
  \def\og{\leavevmode\raise.24ex\hbox{$\scriptscriptstyle\langle\!\langle\>$}}    % guillemets ouvrants
  \def\fg{\leavevmode\raise.24ex\hbox{$\scriptscriptstyle\>\rangle\!\rangle$}}    % guillemets fermants

  \def\d{\,{\rm d}}
  \def\dd{{\rm d}}

  \def\r{{\bb R}}
  \def\CC{{\bb C}}

  \def\HH{{\scal H}}

  \def\O{{\scal O}}
  \def\P{{\scaln P}}

  \def\frac#1#2{{#1\over #2}}
  \def\di#1#2{\sct#1\atop{\sct#2}}

  \def\numero{n$^{\rm o}\thinspace$}

             % carr\'e blanc fin de d\'emo

  \def\numero{n$^{\rm o}\thinspace$}

  \def\¤{\S\thinspace}

  \def\¥{$\bullet$ }
  
  %SYMBOLES MATHS
  
  \def\e{{\rm e}}

  \def\epsilon{\varepsilon}

  \def\phi{\varphi}
  \def\theta{\vartheta}
  \def\rho{\varrho}
  \def\dm{{\textstyle{1\over 2}}}
  \def\txt{\textstyle}
  
  \def\sct{\scriptstyle}
  \def\pf{\noi{\it Proof. }}
  \def\nid{\ifnum\typofr=1\par\noindent{\it D\'emonstration. }\else\pf\fi}
  \def\noi{\noindent}
  \def\rem{\ifnum\typofr=1\noi{\it Remarque.}\ \else\noi{\it Remark.}\ \fi}
  \def\rems{\ifnum\typofr=1\noi{\it Remarques.}\ \else\noi{\it Remarks.}\ \fi}

  \def\1{{\bf 1}}
  \def\|{\Vert}

  \def\leq{\leqslant}
  \def\geq{\geqslant}

  \def\ie{{i.e.\ }}
  \def\eg{{e.g.}}
  \newsymbol\subsetneqq 2324
  \newsymbol\subsetneq 2328
  %% operateurs

  \def\fl#1{\left\lfloor #1 \right\rfloor}

  \def\log{\mathop{\rm log}\nolimits}
  \def\ft#1#2{{\txt{#1\over #2}}}

  %% Grande asterisque operateur

  %% Grande asterisque sur la ligne

  %fleches

  %% Lettres plombees 'pauvres'
  \def\pmb#1{\setbox0=\hbox{#1}%
  \kern-.025em\copy0\kern-\wd0\kern.05em\copy0\kern-\wd0\kern-.025em\raise .0433em\box0 }

  % NOTES EN BAS DE PAGE
  
  % Pour que les accents se placent correctement en mode math en corps 8 et 6
  \skewchar\eighti='177 \skewchar\sixi='177
  \skewchar\eightsy='60 \skewchar\sixsy='60
  
  \def\eightpoint{%
  \textfont0=\eightrm\scriptfont0=\sixrm\scriptscriptfont0=\fiverm
  \def\rm{\fam0\eightrm}%
  \textfont1=\eighti\scriptfont1=\sixi
  \scriptscriptfont1=\fivei\def\oldstyle{\fam1\seveni}%
  \textfont2=\eightsy\scriptfont2=\sixsy\scriptscriptfont2=\fivesy
  \textfont3=\eightex\scriptfont3=\sixex
  \textfont\itfam=\eightit
  \def\it{\fam\itfam\eightit}%
  \textfont\slfam=\eightsl
  \def\sl{\fam\slfam\eightsl}%
  \textfont\bbfam=\eightbb \scriptfont\bbfam=\sixbb\scriptscriptfont\bbfam=\fivebb
  \def\bb{\fam\bbfam\eightbb}%
  \textfont\msafam=\eightmsa\scriptfont\msafam=\sixmsa
  \textfont\scalnfam=\eightscaln
  \def\scaln{\fam\scalnfam\eightscaln}
  \textfont\ttfam=\eighttt
  \def\tt{\fam\ttfam\eighttt}%
\textfont\gotfam=\eightgot
  \textfont\bffam=\eightbf\scriptfont\bffam=\sixbf\scriptscriptfont\bffam=\fivebf
  \def\bf{\fam\bffam\eightbf}%
  \ifx\ITAN\oui\else\textfont\cmmibfam=\eightib
       \scriptfont\cmmibfam=\sixib
        \scriptscriptfont\cmmibfam=\fiveib
         \def\ib{\fam\cmmibfam\eightib}
   \fi
  \textfont\pcapfam=\eightpcap
  \def\pcap{\fam\pcapfam\eightpcap}
  \abovedisplayskip=2pt plus2pt minus 2pt
  \belowdisplayskip=2pt plus1pt minus 2pt
  \abovedisplayshortskip= 1pt plus 2pt minus 1pt
  \belowdisplayshortskip= 1pt plus 2pt minus 1pt
  \smallskipamount=2pt plus 1pt minus 2pt
  \medskipamount=3pt plus 2pt minus 2pt
  \bigskipamount=7pt plus 3pt minus 3pt
  \setbox\strutbox=\hbox{\vrule height 5pt depth 2pt width 0pt}%
  \normalbaselineskip=9pt\normalbaselines\rm}

  \def\({\left(}
  \def\){\right)}
  
  \def\footnoterule{\kern -2pt\hrule width 7truecm\kern 2.4pt}
  
  %reference croisee a des numeros de notes
  \def\xnotedef#1{\definexref{#1}{\noexpand\number\footnotenumber}{Note}}%

  %Fin des definitions pour les notes en bas de page
  
  % PARAGRAPHES EN NEUF POINTS
  
  \def\ninepoint{%
  \textfont0=\ninerm\scriptfont0=\sixrm\scriptscriptfont0=\fiverm
  \def\rm{\fam0\ninerm}%
  \textfont1=\ninei\scriptfont1=\sixi
  \scriptscriptfont1=\fivei\def\oldstyle{\fam1\ninei}%
  \textfont2=\ninesy\scriptfont2=\sixsy\scriptscriptfont2=\fivesy
  \textfont3=\nineex\scriptfont3=\sixex
  \textfont\itfam=\nineit
  \def\it{\fam\itfam\nineit}%
  \textfont\slfam=\ninesl
  \def\sl{\fam\slfam\ninesl}%
  \textfont\bbfam=\ninebb\scriptfont\bbfam=\sixbb\scriptscriptfont\bbfam=\fivebb
  \def\bb{\fam\bbfam\ninebb}%
  \textfont\msafam=\ninemsa\scriptfont\msafam=\sixmsa\scriptscriptfont\msafam=\fivemsa
  \textfont\scalnfam=\ninescaln
  \def\scaln{\fam\scalnfam\ninescaln}
  \textfont\ttfam=\ninett
  \def\tt{\fam\ttfam\ninett}%
  \textfont\bffam=\ninebf\scriptfont\bffam=\sixbf\scriptscriptfont\bffam=\fivebf
  \def\bf{\fam\bffam\ninebf}%
  \abovedisplayskip=3pt plus2pt minus 2pt
  \belowdisplayskip=3pt plus1pt minus 2pt
  \abovedisplayshortskip= 2pt plus 2pt minus 1pt
  \belowdisplayshortskip= 2pt plus 2pt minus 1pt
  \smallskipamount=2pt plus 1pt minus 2pt
  \medskipamount=3pt plus 2pt minus 2pt
  \bigskipamount=7pt plus 3pt minus 3pt
  \setbox\strutbox=\hbox{\vrule height 5pt depth 2pt width 0pt}%
  \normalbaselineskip=11pt plus.3pt minus.2pt\normalbaselines\rm}

  \def\sevenpoint{%
  \textfont0=\sevenrm\scriptfont0=\sixrm\scriptscriptfont0=\fiverm
  \def\rm{\fam0\sevenrm}%
  \textfont1=\seveni\scriptfont1=\sixi
  \scriptscriptfont1=\fivei\def\oldstyle{\fam1\seveni}%
  \textfont2=\sevensy\scriptfont2=\sixsy\scriptscriptfont2=\fivesy
  \textfont3=\sevenex\scriptfont3=\fiveex
  \textfont\itfam=\sevenit
  \def\it{\fam\itfam\sevenit}%
  \textfont\slfam=\sevensl
  \def\sl{\fam\slfam\sevensl}%
  \textfont\bbfam=\sevenbb \scriptfont\bbfam=\sixbb\scriptscriptfont\bbfam=\fivebb
  \def\bb{\fam\bbfam\sevenbb}%
  \textfont\msafam=\sevenmsa\scriptfont\msafam=\sixmsa
  \textfont\scalnfam=\sevenscaln
  \def\scaln{\fam\scalnfam\sevenscaln}
  \textfont\bffam=\sevenbf\scriptfont\bffam=\sixbf\scriptscriptfont\bffam=\fivebf
  \def\bf{\fam\bffam\sevenbf}%
  \textfont\ttfam=\seventt
  \abovedisplayskip=2pt plus2pt minus 2pt
  \belowdisplayskip=2pt plus1pt minus 2pt
  \abovedisplayshortskip= 1pt plus 2pt minus 1pt
  \belowdisplayshortskip= 1pt plus 2pt minus 1pt
  \smallskipamount=2pt plus 1pt minus 2pt
  \medskipamount=3pt plus 2pt minus 2pt
  \bigskipamount=7pt plus 3pt minus 3pt
  \setbox\strutbox=\hbox{\vrule height 5pt depth 2pt width 0pt}%
  \normalbaselineskip=9pt\normalbaselines\rm}

 \def\onzepts{%
 \textfont0=\elevenrm\scriptfont0=\ninerm
 \textfont1=\eleveni\scriptfont1=\ninei
}

\def\douzepts{%
  \textfont0=\twelverm\scriptfont0=\tenrm\def\rm{\fam0\twelverm}%
  \textfont1=\twelvei\scriptfont1=\teni%
  \textfont2=\twelvesy\scriptfont2=\tensy\scriptscriptfont2=\eightsy
  \textfont3=\twelveex
  \textfont\itfam=\twelveti
  \def\it{\fam\itfam\twelveti}%
  \textfont\bffam=\twelvebf\scriptfont\bffam=\tenbf\scriptscriptfont\bffam=\eightbf
  \def\bf{\fam\bffam\twelvebf}%
  \textfont\bbfam=\twelvebb \scriptfont\bbfam=\tenbb
  \def\bb{\fam\bbfam\twelvebb}%
  \textfont\msafam=\twelvemsa\scriptfont\msafam=\tenmsa
  \textfont\scalnfam=\twelvescaln
  \normalbaselineskip=15pt\normalbaselines\rm}

\def\quatorzepts{%
  \textfont0=\fourteenrm\scriptfont0=\twelverm\def\rm{\fam0\fourteenrm}%
  \textfont1=\fourteenib\scriptfont1=\twelveib%
  \textfont2=\fourteensy\scriptfont2=\twelvesy\scriptscriptfont2=\tensy
  \textfont3=\fourteenex
  \textfont\itfam=\fourteenti
  \def\it{\fam\itfam\fourteenti}%
  \textfont\bffam=\fourteenbf\scriptfont\bffam=\twelvebf\scriptscriptfont\bffam=\tenbf
  \def\bf{\fam\bffam\fourteenbf}%
  \textfont\bbfam=\fourteenbb \scriptfont\bbfam=\twelvebb
  \def\bb{\fam\bbfam\fourteenbb}%
  \textfont\msafam=\fourteenmsa\scriptfont\msafam=\twelvemsa
  \textfont\scalnfam=\twelvescaln
  \normalbaselineskip=18pt\normalbaselines\rm}

% Bibliographies, journaux

\def\AA{{\it Acta Arith.}}

%Insertion de figures et illustrations
\def\picture #1 by #2 (#3){\leavevmode\vbox to #2{
     \hrule width #1 height 0pt depth 0pt
      \vfill
       \special{picture #3}}}

\def\illustration #1 by #2 (#3) scaled #4{\dimen1=#2
  \divide\dimen1 by 1000
  \multiply\dimen1 by #4
  \vtop to \dimen1{\dimen1=#1
  \divide\dimen1 by 1000
  \multiply\dimen1 by #4
  \hsize=\dimen1\vss
  \noindent\special{illustration #3 scaled #4}}}

\ifx\couleurs\oui

\fi

\anglais

\optionparag=1
\def\paradouze{oui}
\vsize=255truemm
\voffset=-3truemm
\ifx\optionkeymacros\undefined\else \fi

\catcode`\Œ=\active\defŒ{{\aa}}       % option a
\catcode`\º=\active\defº{\int}        % option b (math mode) 
\catcode`\=\active\def{\c c}        % option c
\catcode`\¶=\active\def¶{\partial}    % option d (math mode)
\catcode`\Ä=\active\defÄ{\oint}       % option f (math mode) ?
\catcode`\Æ=\active\defÆ{\triangle}   % option j (math mode)
\catcode`\Â=\active\defÂ{\neg}        % option l (math mode)
\catcode`\µ=\active\defµ{\mu}         % option m (math mode)
\catcode`\¿=\active\def¿{{\o}}        % option o
\catcode`\¹=\active\def¹{\pi}         % option p (math mode w/ arg.)
\catcode`\Ï=\active\defÏ{{\oe}}       % option q 
\catcode`\§=\active\def§{{\ss}}       % option s 
\catcode`\ =\active\def {\dagger}     % option t  (math mode)
\catcode`\Ã=\active\defÃ{\sqrt}       % option v (math mode w/ arg.)
\catcode`\·=\active\def·{\Sigma}      % option w (math mode)
\catcode`\Å=\active\defÅ{\approx}     % option x (math mode)
\catcode`\½=\active\def½{\Omega}      % option z (math mode)
\catcode`\£=\active\def£{{\it\$}}     % option 3 ($ from italic font)
\catcode`\°=\active\def°{\infty}      % option 5 (math mode)
\catcode`\¤=\active\def¤{{\S}}        % option 6 
\catcode`\¦=\active\def¦{{\P}}        % option 7
\catcode`\¥=\active\def¥{\bullet}     % option 8 
\catcode`\»=\active\def»{\leavevmode\raise.585ex\hbox{\b a}}      % option 9
\catcode`\¼=\active\def¼{\leavevmode\raise.6ex\hbox{\b o}}        % option 0
\catcode`\­=\active\def­{\not=}       % option = (math mode)
\catcode`\²=\active\def²{\leq}        % option , (math mode)
\catcode`\³=\active\def³{\geq}        % option . (math mode)
\catcode`\Ö=\active\defÖ{\div}        % option / (math mode)
\catcode`\É=\active\defÉ{{\dots}}     % option ; 
\catcode`\¾=\active\def¾{{\ae}}       % option '
\catcode`\Ç=\active\defÇ{\og}         % option \ (math mode)
\catcode`\Ò=\active\defÒ{``}          % option [
\catcode`\Á=\active\defÁ{!`}          % option !
\catcode`\¢=\active\def¢{\rlap/c}     % option 4
\catcode`\Ô=\active\defÔ{`}           % option ] 
\catcode`\Õ=\active\defÕ{'}           % shift option ]

% macintosh "shift-option" generated characters

\catcode`\=\active\def{{\AA}}       % shift-option A
\catcode`\'=\active\def'{\c C}        % shift-option C
\catcode`\¯=\active\def¯{{\O}}        % shift-option O
\catcode`\¸=\active\def¸{\Pi}         % shift-option P (math mode)
\catcode`\Î=\active\defÎ{{\OE}}       % shift-option Q
\catcode`\®=\active\def®{{\AE}}       % shift-option '
\catcode`\×=\active\def×{\diamond}    % shift-option V (math mode)
\catcode`\¡=\active\def¡{\accent'27}  % shift-option 8
\catcode`\Ó=\active\defÓ{''}          % shift-option [
\catcode`\±=\active\def±{\pm}         % shift-option = (math mode)
\catcode`\È=\active\defÈ{\fg}         % shift-option \ (math mode)
\catcode`\À=\active\defÀ{?`}          % shift-option / 
\catcode`\Ð=\active\defÐ{--}          % option - (en-dash)
\catcode`\Ñ=\active\defÑ{---}         % shift-option - (em-dash)

% define the macintosh "composite" characters

\catcode`\Š=\active\defŠ{\"a}        % option u, then  a
\catcode`\'=\active\def'{\"e}        % option u, then  e
\catcode`\•=\active\def•{\"{\i}}     % option u, then  i
\catcode`\š=\active\defš{\"o}        % option u, then  o
\catcode`\Ÿ=\active\defŸ{\"u}        % option u, then  u
\catcode`\Ø=\active\defØ{\"y}        % option u, then  y
\catcode`\å=\active\defå{\^A}        %  ^, then  A
\catcode`\€=\active\def€{\"A}        % option u, then  A
\catcode`\…=\active\def…{\"O}        % option u, then  O
\catcode`\†=\active\def†{\"U}        % option u, then  U
\catcode`\‡=\active\def‡{\'a}        % option e, then  a
\catcode`\Ž=\active\defŽ{\'e}        % option e, then  e
\catcode`\'=\active\def'{\'{\i}}     % option e, then  i
\catcode`\—=\active\def—{\'o}        % option e, then  o
\catcode`\œ=\active\defœ{\'u}        % option e, then  u
\catcode`\ƒ=\active\defƒ{\'E}        % option e, then  E
\catcode`\æ=\active\defæ{\^E}        %  ^, then  E
\catcode`\é=\active\defé{\`E}        %  
\catcode`\ˆ=\active\defˆ{\`a}        % option `, then  a
\catcode`\=\active\def{\`e}        % option `, then  e
\catcode`\"=\active\def"{\`{\i}}     % option `, then  i
\catcode`\˜=\active\def˜{\`o}        % option `, then  o
\catcode`\=\active\def{\`u}        % option `, then  u
\catcode`\Ë=\active\defË{\`A}        % option `, then  A
\catcode`\‹=\active\def‹{\~a}        % option n, then  a
\catcode`\–=\active\def–{\~n}        % option n, then  n
\catcode`\›=\active\def›{\~o}        % option n, then  o
\catcode`\Ì=\active\defÌ{\~A}        % option n, then  A
\catcode`\"=\active\def"{\~N}        % option n, then  N
\catcode`\Í=\active\defÍ{\~O}        % option n, then  O
\catcode`\‰=\active\def‰{\^a}        % option i, then  a
\catcode`\=\active\def{\^e}        % option i, then  e
\catcode`\"=\active\def"{\^{\i}}     % option i, then  i
\catcode`\™=\active\def™{\^o}        % option i, then  o
\catcode`\ž=\active\defž{\^u}        % option i, then  u

\let\optionkeymacros\null

\font\tengp=cmmib10

\font\ninegp=cmmib10 at 9pt
               \font\eightgp=cmmib8
              \font\sevengp=cmmib7
%  \font\sixgp=cmmib6
               
               \newfam\gpfam
               \textfont\gpfam=\tengp\scriptfont\gpfam=\sevengp

\font\tenib=cmmib10
         \font\nineib=cmmib10 scaled 900
\font\sevenib=cmmib10 scaled 700
         \font\fiveib=cmmib10 scaled 500
\def\itg{\ib}
\def\0{{\bf 0}}

\def\b{{\itg b}}

\def\FF{{\scal F}}

\def\HH{{\scal H}}

\font\tenib=cmmib10
             \font\nineib=cmmib10 scaled 900
\font\sevenib=cmmib10 scaled 700
             \font\fiveib=cmmib10 scaled 500
\def\itg{\ib}

\def\gR{{\got R}}

\def\auteur{GŽrald Tenenbaum}
 
\def\titrart{On partial derivatives of some summatory functions} 

\hautspages{G. Tenenbaum}{\titrart}
\dateheure
\titrecentre{\titrart}
\bigskip\medskip
\centerline {\auteur}
\bigskip
{\leftskip90mm\it \obeylines
To Helmut Maier,
as a faithful token 
of a long-term companionship.
\par } 
\bigskip
{\eightpoint\leftskip1cm\rightskip1cm
\noi{\bf Abstract.} Let $f$ be a real  arithmetic function and let $g:[1,\infty[\to\r$ be a  smooth function. We describe two emblematic instances in which saddle-point estimates may be used to evaluate the frequency, on the set of integers $n\leqslant x$,  of the event $\{f(n)\leqslant g(n)\}$ from those relevant to the event $\{f(n)\leqslant y\}$. The first example revisits Dickman's historical contribution to the theory of friable integers. The second is concerned with the distribution of the squarefree kernel of an integer.
 \PAR
\medskip
\noi
{\bf Keywords:} friable integers, largest prime factor, squarefree kernel, saddle-point method, local behaviour, semi-asymptotic formulae.\par
\smallskip 
\noi \bf 2020 Mathematics Subject Classification: \rm primary   11N25, 11N37, 11K65; secondary 11N35, 11N64.\par }
\bigskip
\medskip
\paraun{Introduction  and statements of results}
The study of the distribution of an arithmetic function $f$ naturally leads to estimating the  summatory function in two variables
$$\FF(x,y):=\sum_{\di{n\leqslant x}{f(n)\leqslant y}}1.$$ 
However,  passing from an asymptotic formula for $\FF(x,y)$ to an estimate for
$$\HH(x;g):=\sum_{\di{n\leqslant x}{f(n)\leqslant g(n)}}1,\eqdef{Thxg}$$
 may not be straightforward, even when $g$ is a smooth function.
Formally, this task amounts to integrating the partial derivative $\partial\FF(x,g(x))/\partial x$, but,  due to inevitable remainder terms in estimates for $\FF(x,y)$, this is usually out of reach as it stands.
\par 
An obvious path consists in approximating the partial derivative through a discretisation process. This may however turn out to be quite delicate according to available knowledge on the local behaviour of $\FF(x,y)$ with respect to the first variable. A favourable case occurs when $\FF(x,y)$ has been evaluated by the saddle-point method, which  generically provides estimates for the local behaviour---what Erd\H os used to call ``semi-asymptotic formulae''. This note is devoted to describe two instances of this situation. We shall see that this approach  greatly simplifies the computations.
\medskip
The first problem under consideration is related to friable integers, \ie integers free of large prime factors.  Initiated by Dickman  \citer{Di30} in 1930,  the theory of friable integers gradually gained a privileged place in the development of analytic and probabilistic number theory. A direct link with the Riemann hypothesis is made explicit in  \citer{Hi84}. The reader will find a non-exhaustive list of motivations and a description of many fundamental results in the paper \citer{HT86}, in the surveys \citer{HT93}, \citer{Po95}, \citer{Gr08}, in the conference proceedings \citer{DKGL}, and in the monograph \citer{GT15}.  \par 
Let $P^+(n)$ denote the largest prime factor of an integer $n>1$ and put $P^+(1):=1$. Sharp estimates for $$\Psi(x,y):=\sum_{\di{n\leqslant x}{P^+(n)\leqslant y}}1$$
are available in the aforementioned references. However, to  the author's knowledge, none of these revisits the original problem considered by Dickman, that is estimating
$$D(x,u):=\sum_{\di{n\leqslant x}{P^+(n)\leqslant n^{1/u}}}1\quad(x\geqslant 1,\,u\geqslant 1).$$
This question is of type \eqref{Thxg}.\par \goodbreak
 Recall that Dickman's function $\varrho$ is defined as the continuous solution on $\r^+$ to the delay-differential system
$$\normalbaselineskip=15pt\cases{\varrho(v)=1& $(0\leqslant v\leqslant 1)$,\cr
v\varrho'(v)+\varrho(v-1)=0&$(v>1)$.\cr}$$
For $b>0$, $c\geqslant 0$, let $H_{b,c}$ denote the range
$$x\geqslant 3,\quad\e^{(\log_2x)^{b}}<y\leqslant x/(\log x)^{c}.\leqno{(H_{b,c})}$$
Here and throughout, $\log _k$ denotes the $k$th iterated logarithm.
\par \goodbreak
Dickman proved that $D(x,u)\sim x\varrho(u)$ for fixed $u\geqslant 1$ and $x\to\infty$. The best range in which the formula $\Psi(x,y)\sim x\varrho(u)$ ($u:=(\log x)/\log y$) is known to hold is $H_{b,0}$ for any $b>5/3$: this is due to Hildebrand \citer{Hi86}. Thus, the two quantities agree to the first order. We evaluate their difference.
\par 
 Our  statement involves the logarithmic derivative  $$r(v):=-\varrho'(v)/\varrho(v)\quad(v>0).\eqdef{defr}$$
Writing systematically $u:=(\log x)/\log y$, we also use the notation
 $$\gR=\gR(x,y):={\sqrt{\log 2u}\over \sqrt{u}(\log y)^{3/2}}\quad(x\geqslant y\geqslant 2).\eqdef{gR}$$
\Propt{thPsi}{Let $b>5/3$, $c>10$. Uniformly for $(x,y)\in H_{b,c}$,  we have
$$D(x,u)=\Psi(x,y)\Big\{1-{r(u)\over \log y}+O\big(\gR\big)\Big\}.\eqdef{DPsi}$$}
\rems (i) We did not aim at finding the best admissible lower bound for the constant $c$.
\par 
(ii) From Saias' theorem---see \citer{Sa89}, \citeplus{GT15}{ch. III.5}---we have in the same range~$H_{b,c}$
$$\Psi(x,y)=x\varrho(u)+{(\gamma-1)x\varrho'(u)\over \log y}+O\Big({x\varrho(u)(\log 2u)^2\over (\log y)^2}\Big),$$
where $\gamma$ denotes Euler's constant. We hence deduce from \eqref{DPsi} that, still in $H_{b,c}$,
$$D(x,u)=x\varrho(u)+{\gamma x\varrho'(u)\over \log y}+O\big(x\varrho(u)\gR_1\big),\eqdef{Drho}$$
with $\gR_1:=\gR+(\log 2u)^2/(\log y)^2\ll(\log 2u)^{7/6}/(\log y)^{3/2}.$
\medskip
(iii) It is easy to deduce from Buchstab's identity that, for $\sqrt{x}<y\leqslant x$, we have
$$D(x,u)=\Psi(x,y)-\sum_{p\leqslant y}\big\{p^{u-1}+O(1)\big\}=x\varrho(u)+{\gamma x\varrho'(u)\over \log y}+O\bigg({y\over \log x}+{x\over (\log x)^2}\bigg).$$
\medskip
The following statement is a straightforward consequence of \ref{thPsi}.
\Propc{corPsi}{We have
$$\sum_{1<n\leqslant x}{\log n\over \log P^+(n)}=\e^\gamma x-{\gamma\e^\gamma x\over \log x}+O\Big({x\over (\log x)^{3/2}}\Big)\qquad (x\geqslant 3).\eqdef{ln/lP}$$}
If $\log n$ is replaced by $\log x$ in the left-hand side of \eqref{ln/lP}, the coefficient $-\gamma\e^\gamma$ changes to~$(1-\gamma)\e^\gamma$.\bigskip
Our second investigation is about the distribution of the squarefree kernel of an integer,~\ie
$$k(n):=\prod_{p|n}p\quad(n\geqslant 1),$$
where, here and in the sequel, $p$ denotes a prime number. In a recent article, BrŸdern \& Robert~\citer{BR24} obtained an estimate for
$$S(x;\vartheta,\alpha):=\sum_{\di{n\leqslant x}{k(n)\leqslant n^{\vartheta}(\log n)^\alpha}}1$$
for fixed $\vartheta\in]0,1[$, $\alpha\in\r$. The following statement significantly improves on their result regarding uniformity ($\vartheta$ and $\alpha$ will not assumed to be fixed) and domain of validity ($\vartheta$ will be allowed to approach 0 or 1). Moreover, replacing $n^{\vartheta}(\log n)^\alpha$ by another smooth function would easily follow by the same approach. The main ingredients of the proof, which turns out to be very short, are estimates obtained in \citer{RT13} for the local behaviour of the summatory function
$$N(x,y):=\sum_{\di{n\leqslant x}{k(n)\leqslant y}}1.$$
 These estimates are all established by the saddle-point method. \par \goodbreak
We need the following notation.
Let $\psi(n):=\prod_{p|n}(p+1)$ $(n\geqslant 1)$, and define
$$F(t):={6\over \pi^2}\sum_{n\geqslant 1}{\min(1,\e^t/n)\over n\psi(n)}\qquad (t>0).$$
Furthermore, consider $g(\sigma):=\sum_{p}\log \{1+(1-p^{\sigma-1})/p(p^\sigma-1)\}$ and, for $t\geqslant 1$, let $\sigma_t$ denote the solution to the equation $g'(\sigma)+t=0$. From \citeplus{RT13}{(2.9)} we have 
$$
\sigma_t=
\sqrt{\frac{2}{t\log t}}\bigg\{1+\sum_{1\leqslant
k\leqslant K}{P_k(\log_2t)\over (\log
t)^k}+O_K\bigg({(\log_2 t)^{K+1}\over (\log
t)^{K+1}}\bigg)\bigg\}\quad
(t\to\infty),\eqdef{fa-sigv}$$
where $P_k$ is a polynomial of degree at most $k$.
In particular, $$P_1(z)=\dm(z-\log 2), \quad
P_2(z)=\ft38z^2-(\ft34\log 2+\dm)z+\dm\log
2+\ft38(\log 2)^2+\ft23\pi^2.\eqdef{pols-dasig}$$
\goodbreak
\Propt{thN}{Let $A>0$, $c\in]0,\dm[$. Uniformly for $x\geqslant 3$, $1/(\log x)^c\leqslant \vartheta\leqslant 1-1/(\log x)^c$, $|\alpha|\leqslant A$, $y=x^\vartheta(\log x)^\alpha$, $v:=\log (x/y)$, we have
$$\eqalign{S(x;\vartheta,\alpha)&={yF(v)\sigma_v\over \vartheta}\bigg\{1-{\sigma_v\over \vartheta}+O\bigg(\frac{\sigma_v^2}{\vartheta^2}+\sqrt{\log v\over v}\bigg)\bigg\}\cr&={N(x,y)\sigma_v\over \vartheta}\bigg\{1-{\sigma_v\over \vartheta}+O\bigg(\frac{\sigma_v^2}{\vartheta^2}+\sqrt{\log v\over v}\bigg)\bigg\}.}\eqdef{estSxth}$$}
Since $\sigma_v/\vartheta\ll(\log x)^{c-1/2}$, we see that carrying~\eqref{fa-sigv} back into~\eqref{estSxth} yields an asymptotic expansion of the left-hand side with general term $c_{\ell k}(\log_2v)^\ell/(\log v)^k$,  $0\leqslant \ell\leqslant k$,  $c_{\ell k}\in\r$.
\bigskip\bigskip
\paraun{Proof of \ref{thPsi} }
Let us assume throughout that $x$ is arbitrarily large since the stated estimates are otherwise trivial. Recall the notation $u:=(\log x)/\log y$ and the definition of the function $r$ in \eqref{defr}. Put $Z(s):=(s-1)\zeta(s)/s$ $(s\in\CC)$, where $\zeta(s)$ is the Riemann zeta function. From a slight modification of \citeplus{BT05}{prop. 4.1(i)}, we derive that, for $(x,y)\in(H_{b,c})$, $\beta:=1-r(u)/\log y$, $b>5/3$, $c>8$, we have
$$\Psi(x,y)=x\varrho(u)Z(\beta)\Big\{1+O\Big({u\over (\log x)^2}\Big)\Big\}.\eqdef{main}$$
The  alteration mentioned above solely concerns the lower bound on $c$, which is anyway not optimal. It stems from the fact that no coprimality is required here, which slightly simplifies the argument.  
\par 
 Estimate \eqref{main}, which is proved through a saddle-point analysis, is the main tool in the proof of~\eqref{DPsi}, inasmuch it provides an estimate to the second order.
 \par 
 As mentioned in the introduction, we mimic partial derivatives $\partial\Psi(x,y)/\partial x$ at $y=x^{1/u}$. Thus, given a small parameter $\varepsilon=\varepsilon_x>0$ to be chosen later, we put $x_k:=x\e^{-k\varepsilon}$, $y_k:=x_k^{1/u}$ $(k\geqslant 0)$ and write, for some integer $K\geqslant 1$,
$$\leqalignno{D^-(x,u)&:=\sum_{k<K}\big\{\Psi\big(x_k,y_{k+1}\big)-\Psi\big(x_{k+1},y_{k+1}\big)\big\},&\cr
D^+(x,u)&:= \sum_{k<K}\big\{\Psi\big(x_k,y_k\big)-\Psi\big(x_{k+1},y_k\big)\big\}+\Psi\big(x_K,y_K\big),&\cr
&D^-(x,u)\leqslant D(x,u)\leqslant D^+(x,u).&\eqdef{encD} \cr}$$
\par 
Select $K:=\fl{2(\log_2x)/\varepsilon}$, so that $\Psi\big(x_K,y_K\big)\ll x\varrho(u)/(\log x)^2$, a quantity that may be absorbed by the remainder term of \eqref{DPsi}. This choice implies in particular that
 $\log x_k\asymp \log x$ $(0\leqslant k\leqslant K).$
\par 
Observe that $(x,y)\in H_{b,c}$ with $b>5/3$, $c>10$ implies $(x_k,y_k)\in H_{b',c'}$ for some $b'>5/3$, $c'>8$. Hence our hypotheses enable applying \eqref{main} for $(x_k,y_k)$, viz.
$$\Psi\big(x_k,y_k\big)=x_kZ(\beta_k)\varrho(u)\Big\{1+O\Big({1\over u(\log y)^2}\Big)\Big\}\quad(0\leqslant k<K),\eqdef{Pk}$$
with $\beta_k:=1-r(u)/\log y_k.$ By the same device, we have
$$\Psi\big(x_{k+1},y_{k}\big)=x_{k+1}Z(\gamma_{k})\varrho(v_k)\Big\{1+O\Big({1\over u(\log y)^2}\Big)\Big\}\quad(0\leqslant k<K),\eqdef{Pk+1}$$
with $$v_k:={\log x_{k+1}\over \log y_k}=u-{\varepsilon\over \log y_k},\quad\gamma_k:=1-{r(v_k)\over\log y_k}\cdot$$
\par 
If $\xi(v)$ denotes the solution to the equation $\e^\xi=1+v\xi$ for $v\neq1$ and $\xi(1):=0$, we have \citer{BT05}, for $v\geqslant 3$,
$$\eqalign{&\xi(v)=\log (v\log v)+O\Big({\log_2v\over \log v}\Big),\quad\xi'(v)={1\over v}+{1\over v\log v}+O\Big({\log_2v\over v(\log v)^2}\Big),\cr
&r(v)-\xi(v)\ll{1\over v},\quad r'(v)-\xi'(v)\ll{1\over v^2}.\cr}\eqdef{xir}$$
Therefore $r(v_k)-r(u)\ll\varepsilon/\log x$, and
$\gamma_k-\beta_k\ll\varepsilon / \{u(\log y)^2\}.$
Applying the  Taylor-Lagrange formula \citeplus{GT15}{(5.115)} for $\varrho(v)$ at order 1, we get 
$$\varrho(v_k)=\varrho(u)-{\varepsilon \varrho'(u)\over \log y_k}+O\Big({\varepsilon^2\varrho(u)(\log 2u)^2\over (\log y)^2}\Big),\quad Z(\gamma_k)=Z(\beta_k)+O\Big({\varepsilon \over u(\log y)^2}\Big),$$
under the extra condition  $u>1+2\varepsilon/\log y_k$, which is actually  implied by the hypothesis $H_{b,c}$.\par 
Carrying back into \eqref{Pk+1} and assuming  $\varepsilon\ll1/\{\sqrt{u}\log 2u\}$, we derive in turn
$$\eqalign{\Psi\big(x_{k+1},y_{k}\big)&=x_{k+1}Z(\beta_k)\varrho(u)\bigg\{1+{\varepsilon r(u)\over \log y_k}+O\Big({ 1\over u(\log y)^2}\Big)\bigg\},\cr
\Psi\big(x_k,y_k\big)-\Psi\big(x_{k+1},y_k\big)&=x_k(1-\e^{-\varepsilon})Z(\beta_k)\varrho(u)+{\varepsilon \e^{-\varepsilon}x_kZ(\beta_k)\varrho'(u)\over \log y_k}+O\Big({x_k\varrho(u)\over u(\log y)^2}\Big).\cr}$$
Now 
$$\log y_k=\log y-k\varepsilon/u,\quad Z(\beta_k)=Z(\beta)+O\Big({k\varepsilon\log 2u\over u(\log y)^2}\Big), $$
whence 
$$\eqalign{\Psi\big(x_k,y_k\big)-&\Psi\big(x_{k+1},y_k\big)\cr
&=x_k(1-\e^{-\varepsilon})Z(\beta)\varrho(u)+{\varepsilon x_kZ(\beta)\e^{-\varepsilon}\varrho'(u)\over \log y}+O\bigg({x_k\varrho(u)\{1+k\varepsilon^2(\log 2u)^2\}\over u(\log y)^2}\bigg).\cr}$$
Summing over $k$ yields 
$$D^+(x,u)= {x\varrho(u)Z(\beta)}+{Z(\beta)x\varrho'(u)\over \log y}+O\bigg({\varepsilon x\varrho(u)\log (2u)\over \log y}+{x\varrho(u)\over \varepsilon u(\log y)^2}+{x\varrho(u)(\log 2u)^2\over u(\log y)^2}\bigg).$$
For the quasi-optimal choice $\varepsilon:=1/\sqrt{(\log x)\log 2u}$, we get, with notation \eqref{gR},
$$\eqalign{D^+(x,u)&= {x\varrho(u)Z(\beta)}+{Z(\beta)x\varrho'(u)\over \log y}+O\big(x\varrho(u)\gR\big)=x\varrho(u)Z(\beta)\Big\{1-{r(u)\over \log y}+O\big(\gR\big)\Big\},\cr}$$
which is compatible with \eqref{DPsi}.
Parallel computations provide the same estimate for $D^-(x,u)$.
\medskip
\bigskip
\paraun{Proof of  \ref{corPsi}}
Let $D(x)$ denote the left-hand side of \eqref{ln/lP}. We plainly have
$$D(x)=\int_0^\infty \{D(x,u)-1\}\d u.\eqdef{intD}$$
From the bound $1\leqslant D(x,u)\leqslant \Psi\big(x,x^{1/u}\big)\ll x\e^{-u/2}$ $(u\geqslant 1,\,x\geqslant 1)$ \citeplus{GT15}{th. III.5.1}, we see that the contribution of the range $u>4\log_2x$ may be encompassed by the stated remainder. 
\par 
Let $\varepsilon_x:=11(\log_2x)/\log x$. For large $x$ and $1+\varepsilon_x<u\leqslant 4\log_2x$, formula \eqref{Drho} is applicable. It is also trivially valid for $0\leqslant u\leqslant 1$ at the cost of replacing the error term by $O(1)$. When $1<u\leqslant 1+\varepsilon_x$, we have, writing $X:=\{x/\log x\}^{1/(1+\varepsilon_x)}$,
$$0\leqslant x-D(x,u)=x\varrho(u)-D(x,u)\leqslant \sum_{X<p\leqslant x}\fl{x\over p}+{x\over \log x}\ll{x\log_2x\over \log x}\cdot$$
Therefore
$$\int_0^{1+\varepsilon_x}\Big\{D(x,u)-1-x\varrho(u)-{\gamma xu\varrho'(u)\over \log x}\Big\}\d u\ll {x(\log_2x)^2\over (\log x)^2}\cdot$$
Applying \eqref{Drho} to evaluate the contribution  of the range $1+\varepsilon_x\leqslant u\leqslant 4\log_2x$ to the integral in~\eqref{intD}, we get
$$\int_0^{4\log_2x}\Big\{D(x,u)-1-x\varrho(u)+{\gamma x\varrho(u-1)\over \log x}\Big\}\d u\ll {x\over (\log x)^{3/2}}\cdot$$
Now \eqref{ln/lP} follows by extending the last integrals to infinity and appealing to the formula $\int_0^\infty\varrho(u)\d u=\e^\gamma$---see \eg\ \citeplus{GT15}{(III.5.45)}.
\medskip\bigskip
\paraun{Proof of \ref{thN}}
As previously we assume $x$ arbitrarily large and mimic partial derivatives $\partial N(x,y)/\partial x$ at $y:=x^\vartheta(\log x)^\alpha$. 
Our main tool will be the following estimate \citeplus{RT13}{(3.25)}, valid for fixed $b>\dm$ and uniformly for $x\geqslant 3$,  $\e^{(\log x)^b}<y\leqslant x$, $v:=\log (x/y)$, $\eta_x:=\sqrt{2/(\log x)\log_2x}$:
$$N(x,y)=yF(v)\big\{1+O\big(y^{-\eta_x}\big)\big\}.\eqdef{fondN}$$
Given $\varepsilon=\varepsilon_x\in]0,1[$ to be determined later, we put $x_k:=x\e^{-k\varepsilon}$, $y_k:=x_k^\vartheta(\log x_k)^\alpha$ $(k\geqslant 0)$, and, given an integer $K\in[1, (\log x)/2\varepsilon]$, write
$$\leqalignno{S^-(x;\vartheta,\alpha)&:=\sum_{k<K}\big\{N\big(x_k,y_{k+1}\big)-N\big(x_{k+1},y_{k+1}\big)\big\},&\cr
S^+(x;\vartheta,\alpha)&:= \sum_{k<K}\big\{N\big(x_k,y_k\big)-N\big(x_{k+1},y_k\big)\big\}+N\big(x_K,y_K\big),&\cr
&S^-(x;\vartheta,\alpha)\leqslant S(x;\vartheta,\alpha)\leqslant S^+(x;\vartheta,\alpha).&\eqdef{encS} \cr}$$
We note right-away that, by \eqref{fondN} and \citeplus{RT13}{(8.13)} in the form
$$F(v+h)\ll\e^{h\sigma_v}F(v)\quad\big(v\geqslant 2,\,h+v\geqslant 0\big),$$
and  since  $\log(x_K/y_K)=v-(1-\vartheta)K\varepsilon+O(1)$,  we have
$$N(x_K,y_K)\ll y_K\e^{-(1-\vartheta)\varepsilon K\sigma_v}F(v)\ll\e^{-\varepsilon \vartheta K}yF(v).$$
Therefore, selecting $K:=\fl{2(\log v)/\varepsilon\vartheta}$  ensures that the last term in the upper bound for $S^+(x;\vartheta,\alpha)$ is absorbable by the error term of \eqref{estSxth}. Note that this implies $v_k:=\log (x_k/y_k)\asymp v$ $(0\leqslant k\leqslant K)$.
\par \goodbreak
Next, we apply \eqref{fondN} to $(x_k,y_k)$ and $(x_{k+1},y_k)$. This yields, with, say, $R:=\e^{-(\log x)^{1/4-c/2}}$,
$$\eqalign{N(x_k,y_k)&=y_kF(v_k)\big\{1+O\big(R\big)\big\},\cr N(x_{k+1},y_k)&=y_kF(v_k-\varepsilon)\big\{1+O\big(R\big)\big\}.\cr}$$
Now, \citeplus{RT13}{(8.16)} furnishes
$$F(w)-F(w-h)=h \sigma_{w}F(w)\bigg\{1+O\bigg({|h|+\log w\over \sqrt{w\log w}}\bigg)\bigg\}\qquad \big(h\ll\sqrt{w\log w}\big),\eqdef{accrF}$$
from which, specializing $w:=v_k$, $h:=\varepsilon$, we derive
$$N(x_k,y_k)-N(x_{k+1},y_k)=y_k\varepsilon\sigma_{v_k}F(v_k)\bigg\{1+O\bigg(\sqrt{\log v\over v}+{R\over \varepsilon\sigma_v}\bigg)\bigg\}.$$
The penultimate step consists in substituting $v$ to $v_k$ in the above formula. By \citeplus{RT13}{(8.4)}, we know that $\dd\sigma_v/\dd v\ll v^{-3/2}(\log v)^{-1/2}$, hence $$\sigma_{v_k}-\sigma_v\ll{k\varepsilon\over v\sqrt{v\log v}}\ll{k\varepsilon\sigma_v\over v},\quad F(v_k)=F(v)\bigg\{1-k\varepsilon\sigma_v+O\bigg(k^2\varepsilon^2\sigma_v^2+\sqrt{\log v \over v}\bigg)\bigg\},$$
where the last estimate follows from \eqref{accrF} with $w:=v$, $h:=(1-\vartheta)\varepsilon k-\alpha\log (1-\varepsilon k/\log x)$, observing that $\vartheta\varepsilon k\sigma_v\ll\sqrt{(\log v)/v}$.
Finally, we arrive at
$$N(x_k,y_k)-N(x_{k+1},y_k)=y_k\varepsilon\sigma_{v}F(v)\bigg\{1-k\varepsilon\sigma_v+O\bigg(k^2\varepsilon^2\sigma_v^2+\sqrt{\log v\over v}+{R\over \varepsilon \sigma_v}\bigg)\bigg\}.$$ 
Selecting $\varepsilon=\sqrt{(\log v)/v}$ and summing over $k$ provides the required estimate for $S^+(x;\vartheta,\alpha)$. The corresponding formula for $S^-(x;\vartheta,\alpha)$ is proved by parallel computations.
\bigskip
\noi{\bf Acknowledgement.} The author takes pleasure in thanking RŽgis de la Bretche for a fruitful conversation on the topic of this article, and the referee for pertinent suggestions.

\bigskip\bigskip
 \centerline{\twelvebf References}
\bigskip
 {\eightpoint\leftskip9mm\rightskip5mm
\bibtem{BT05} R. de la Bretche \& G. Tenenbaum, PropriŽtŽs statistiques des entiers friables, {\it Ramanujan J. \bf9} (2005),
139--202.\par  
\bibtem{BR24} J. BrŸdern \& O. Robert, On the distribution of powered numbers, {\it Open Math. \bf 22} (2024) no.\thinspace1, 20240007.\par 
\bibtem{DKGL} J.-M. De Koninck, A. Granville \& F. Luca (eds.) {\it Anatomy of integers}, Papers from the CRM Workshop held at the UniversitŽ de MontrŽal, MontrŽal, QC, March 13Ð17, 2006, CRM Proceedings \& Lecture Notes, 46, American Mathematical Society, Providence, RI, 2008. viii+297~pp.\par 
\bibtem{Di30} K. Dickman,
 On the frequency of numbers containing prime factors of a certain relative
magnitude, {\it Ark. Math. Astr. Fys. \bf 22} (1930), 1--14.
\par 
\bibtem{Gr08}A. Granville, Smooth numbers: computational number theory and beyond, in: {\it Algorithmic number theory: lattices, number fields, curves and cryptography}, Math. Sci. Res. Inst. Publ. 44, Cambridge Univ. Press, Cambridge, 2008.\par 
\bibtem{Hi84}A. Hildebrand, Integers free of large prime factors and the Riemann hypothesis, {\it Mathematika \bf31} (1984), 58-271.
\par 
\bibtem{Hi86}A. Hildebrand, On the number of positive integers $\leqslant  x$ and free of
prime factors $>y$, {\it J. Number Theory   \bf 22} (1986), 289--307.
\par 
\bibtem{HT86}A. Hildebrand \& G. Tenenbaum, On integers free of large prime factors, {\it Trans. Amer.
Math. Soc. \bf 296} (1986), 265--290.
\par 
\bibtem {HT93} A. Hildebrand \& G. Tenenbaum, Integers without large prime factors, {\it J. Th\'eorie des
Nombres de   Bordeaux \bf 5} (1993), 411-484. 
\par \bibtem{Po95} C. Pomerance, The role of smooth numbers in number-theoretic algorithms, in : {\it Proceedings of the International Congress of Mathematicians}, Vol. 1, 2 (ZŸrich, 1994), 411Ð422, BirkhŠuser, Basel, 1995.\par 

\bibtem{RT13} O. Robert \& G. Tenenbaum, Sur la rŽpartition du noyau d'un entier, {\it Indag. Math. \bf 24} (2013), 802--914.\par 
\bibtem{Sa89} ƒ. Saias, Sur le nombre des entiers sans grand facteur premier, {\it J. Number Theory \bf
32} \numero 1 (1989), 78--99.\par 
\bibtem{GT15} G. Tenenbaum, {\it Introduction to analytic and probabilistic number theory}, Graduate Studies in Mathematics 163, Amer. Math. Soc. 2015; see also {\it Introduction ˆ la thŽorie analytique et probabiliste des
nombres}, 5th ed., Dunod, Sciences Sup 2022, xvii+547 pp.

\bigskip\bigskip
{\leftskip2mm\rightskip-2cm\sevenrm
 \obeylines \baselineskip=8pt
 GŽrald Tenenbaum
Institut \'Elie Cartan
Universit\'e de Lorraine
 BP 70239\par
54506 Vand\oe uvre-ls-Nancy Cedex
 France
\smallskip
{\seventt gerald.tenenbaum@univ-lorraine.fr}
\par}

\vfill\eject

\end